\documentclass[12pt]{article}%
\usepackage{graphicx}
\usepackage[intlimits]{amsmath}
\usepackage{amsfonts}
\usepackage{amssymb}%
 \makeatletter 
\newtheorem{theorem}{Theorem}

\newenvironment{proof}[1][Proof]{\noindent{\textbf {#1}  }}  {\hfill$\Box$\bigskip}

\begin{document}

\title{The Energy of Graphs and Matrices}
\author{Vladimir Nikiforov\\Department of Mathematical Sciences, University of Memphis, \\Memphis TN 38152, USA, e-mail:\textit{ vnikifrv@memphis.edu}}
\maketitle

\begin{abstract}
Given a complex $m\times n$ matrix $A,$ we index its singular values as
$\sigma_{1}\left(  A\right)  \geq\sigma_{2}\left(  A\right)  \geq...$ and call
the value $\mathcal{E}\left(  A\right)  =\sigma_{1}\left(  A\right)
+\sigma_{2}\left(  A\right)  +...$ the \emph{energy }of $A,$ thereby extending
the concept of graph energy, introduced by Gutman. Let $2\leq m\leq n,$ $A$ be
an $m\times n$ nonnegative matrix with maximum entry $\alpha$, and $\left\Vert
A\right\Vert _{1}\geq n\alpha$. Extending previous results of Koolen and
Moulton for graphs, we prove that
\[
\mathcal{E}\left(  A\right)  \leq\frac{\left\Vert A\right\Vert _{1}}{\sqrt
{mn}}+\sqrt{\left(  m-1\right)  \left(  tr\left(  AA^{\ast}\right)
-\frac{\left\Vert A\right\Vert _{1}^{2}}{mn}\right)  }\leq\alpha\frac{\sqrt
{n}\left(  m+\sqrt{m}\right)  }{2}.
\]
Furthermore, if $A$ is any nonconstant matrix, then
\[
\mathcal{E}\left(  A\right)  \geq\sigma_{1}\left(  A\right)  +\frac{tr\left(
AA^{\ast}\right)  -\sigma_{1}^{2}\left(  A\right)  }{\sigma_{2}\left(
A\right)  }.
\]

Finally, we note that Wigner's semicircle law implies that%
\[
\mathcal{E}\left(  G\right)  =\left(  \frac{4}{3\pi}+o\left(  1\right)
\right)  n^{3/2}%
\]
for almost all graphs $G$.


\end{abstract}

Our notation is standard (e.g., see \cite{Bol98}, \cite{CDS80}, and
\cite{HoJo88}); in particular, we write $M_{m,n}$ for the set of $m\times n$
matrices with complex entries, and $A^{\ast}$ for the Hermitian adjoint of
$A.$ The singular values $\sigma_{1}\left(  A\right)  \geq\sigma_{2}\left(
A\right)  \geq...$ of a matrix $A$ are the square roots of the eigenvalues of
$AA^{\ast}.$ Note that if $A\in M_{n,n}$ is a Hermitian matrix with
eigenvalues $\mu_{1}\left(  A\right)  \geq...\geq\mu_{n}\left(  A\right)  $,
then the singular values of $A$ are the moduli of $\mu_{i}\left(  A\right)  $
taken in descending order.

For any $A\in M_{m,n},$ call the value $\mathcal{E}\left(  A\right)
=\sigma_{1}\left(  A\right)  +...+\sigma_{n}\left(  A\right)  $ the
\emph{energy }of $A$. Gutman \cite{Gut78} introduced $\mathcal{E}\left(
G\right)  =\mathcal{E}\left(  A\left(  G\right)  \right)  ,$ where $A\left(
G\right)  $ is the adjacency matrix of a graph $G;$ in this narrow sense
$\mathcal{E}\left(  A\right)  $ has been studied extensively (see, e.g.,
\cite{Bal04}, \cite{Gut99}, \cite{KoMo01}, \cite{KoMo03}, \cite{RaTi04},
\cite{Shp06}, and \cite{StSt05}). In particular, Koolen and Moulton
\cite{KoMo01} proved the following sharp inequalities for a graph $G$ of order
$n$ and size $m\geq n/2,$
\begin{equation}
\mathcal{E}\left(  G\right)  \leq2m/n+\sqrt{\left(  n-1\right)  \left(
2m-\left(  2m/n\right)  ^{2}\right)  }\leq\left(  n/2\right)  \left(
1+\sqrt{n}\right)  . \label{KM}%
\end{equation}
Moreover, Koolen and Moulton conjectured that for every $\varepsilon>0$, for
almost all $n\geq1,$ there exists a graph $G$ with $\mathcal{E}\left(
G\right)  \geq\left(  1-\varepsilon\right)  \left(  n/2\right)  \left(
1+\sqrt{n}\right)  $.

In this note we give upper and lower bounds on $\mathcal{E}\left(  A\right)  $
and find the asymptotics of $\mathcal{E}\left(  G\right)  $ of almost all
graphs $G$. We first generalize inequality (\ref{KM}) in the following way.

\begin{theorem}
\label{tupb}If $m\leq n,$ $A$ is an $m\times n$ nonnegative matrix with
maximum entry $\alpha,$ and $\left\Vert A\right\Vert _{1}\geq n\alpha,$ then
\begin{equation}
\mathcal{E}\left(  A\right)  \leq\frac{\left\Vert A\right\Vert _{1}}{\sqrt
{mn}}+\sqrt{\left(  m-1\right)  \left(  tr\left(  AA^{\ast}\right)
-\frac{\left\Vert A\right\Vert _{1}^{2}}{mn}\right)  }. \label{upb}%
\end{equation}

\end{theorem}

From here we derive the following absolute upper bound on $\mathcal{E}\left(
A\right)  .$

\begin{theorem}
\label{taupb}If $m\leq n$ and $A$ is an $m\times n$ nonnegative matrix with
maximum entry $\alpha,$ then,
\begin{equation}
\mathcal{E}\left(  A\right)  \leq\alpha\frac{\left(  m+\sqrt{m}\right)
\sqrt{n}}{2}. \label{aupb}%
\end{equation}

\end{theorem}

Note that Theorems \ref{tupb} and \ref{taupb} improve on the bounds for the
energy of bipartite graphs given in \cite{KoMo03}.

On the other hand, for every $A\in M_{m,n},$ $\left(  m,n\geq2\right)  ,$ we
have $\sigma_{1}^{2}\left(  A\right)  +\sigma_{2}^{2}\left(  A\right)
+...=tr\left(  AA^{\ast}\right)  ,$ and so
\[
tr\left(  AA^{\ast}\right)  -\sigma_{1}^{2}\left(  A\right)  =\sigma_{2}%
^{2}+...+\sigma_{m}^{2}\leq\sigma_{2}\left(  A\right)  \left(  \mathcal{E}%
\left(  A\right)  -\sigma_{1}\left(  A\right)  \right)  .
\]
Thus, if $A$ is a nonconstant matrix, then
\begin{equation}
\mathcal{E}\left(  A\right)  \geq\sigma_{1}\left(  A\right)  +\frac{tr\left(
AA^{\ast}\right)  -\sigma_{1}^{2}\left(  A\right)  }{\sigma_{2}\left(
A\right)  }. \label{lowb}%
\end{equation}

If $A$ is the adjacency matrix of a graph, this inequality is tight up to a
factor of $2$ for almost all graphs. To see this, recall that the adjacency
matrix $A\left(  n,1/2\right)  $ of the random graph $G\left(  n,1/2\right)  $
is a symmetric matrix with zero diagonal, whose entries $a_{ij}$ are
independent random variables with $E\left(  a_{ij}\right)  =1/2,$ $Var\left(
a_{ij}^{2}\right)  =1/4=\sigma^{2},$ and $E\left(  a_{ij}^{2k}\right)
=1/4^{k}$ for all $1\leq i<j\leq n,$ $k\geq1.$ The result of F\"{u}redi and
Koml\'{o}s \cite{FuKo81} implies that, with probability tending to 1,
\begin{align*}
\sigma_{1}\left(  G\left(  n,1/2\right)  \right)   &  =\left(  1/2+o\left(
1\right)  \right)  n,\text{ }\\
\sigma_{2}\left(  G\left(  n,1/2\right)  \right)   &  <\left(  2\sigma
+o\left(  1\right)  \right)  n^{1/2}=\left(  1+o\left(  1\right)  \right)
n^{1/2}.
\end{align*}
Hence, inequalities (\ref{KM}) and (\ref{lowb}) imply that%
\[
\left(  1/2+o\left(  1\right)  \right)  n^{3/2}>\mathcal{E}\left(  G\right)
>\left(  1/2+o\left(  1\right)  \right)  n+\frac{\left(  1/4+o\left(
1\right)  \right)  n^{2}}{\left(  1+o\left(  1\right)  \right)  n^{1/2}%
}=\left(  1/4+o\left(  1\right)  \right)  n^{3/2}%
\]
for almost all graphs $G$.

Moreover, Wigner's semicircle law \cite{Wig58} (we use the form given by
Arnold \cite{Arn67}, p. 263), implies that
\[
\mathcal{E}\left(  A\left(  n,1/2\right)  \right)  n^{-1/2}=n\left(  \frac
{2}{\pi}\int_{-1}^{1}\left\vert x\right\vert \sqrt{1-x^{2}}dx+o\left(
1\right)  \right)  =\left(  \frac{4}{3\pi}+o\left(  1\right)  \right)  n,
\]
and so $\mathcal{E}\left(  G\right)  =\left(  \frac{4}{3\pi}+o\left(
1\right)  \right)  n^{3/2}$ for almost all graphs $G$.

\begin{proof}
[\textbf{Proof of Theorem \ref{tupb}}]We adapt the proof of (\ref{KM}) in
\cite{KoMo01}. Letting $\mathbf{i}$ to be the all ones $m$-vector, Rayleigh's
principle implies that $\sigma_{1}^{2}\left(  A\right)  m\geq\left\langle
AA^{\ast}\mathbf{i},\mathbf{i}\right\rangle ;$ hence, after some algebra,
$\sigma_{1}\left(  A\right)  \geq\left\Vert A\right\Vert _{1}/\sqrt{mn}.$ The
AM-QM inequality implies that,
\[
\mathcal{E}\left(  A\right)  -\sigma_{1}\left(  A\right)  \leq\sqrt{\left(
m-1\right)  \sum_{i=2}^{n}\sigma_{i}^{2}\left(  A\right)  }=\sqrt{\left(
m-1\right)  \left(  tr\left(  AA^{\ast}\right)  -\sigma_{1}^{2}\left(
A\right)  \right)  }.
\]
The function $x\rightarrow x+\sqrt{\left(  m-1\right)  \left(  tr\left(
AA^{\ast}\right)  -x^{2}\right)  }$ is decreasing if $\sqrt{tr\left(
AA^{\ast}\right)  /m}\leq x\leq\sqrt{tr\left(  AA^{\ast}\right)  };$ hence, in
view of
\[
tr\left(  AA^{\ast}\right)  =\sum_{j=1}^{n}\sum_{k=1}^{m}\left\vert
a_{kj}\right\vert ^{2}=\sum_{j=1}^{n}\sum_{k=1}^{m}a_{kj}^{2}\leq\alpha
\sum_{j=1}^{n}\sum_{k=1}^{m}a_{kj}=\alpha\left\Vert A\right\Vert _{1},
\]
we find that $\sqrt{tr\left(  AA^{\ast}\right)  /m}\leq\left\Vert A\right\Vert
_{1}/\sqrt{mn}$, and inequality (\ref{upb}) follows.
\end{proof}

\begin{proof}
[\textbf{Proof of Theorem \ref{taupb}}]If $\left\Vert A\right\Vert _{1}\geq
n\alpha,$ then Theorem \ref{tupb} and $tr\left(  AA^{\ast}\right)  \leq
\alpha\left\Vert A\right\Vert _{1}$ imply that%
\[
\mathcal{E}\left(  A\right)  \leq\frac{\left\Vert A\right\Vert _{1}}{\sqrt
{mn}}+\sqrt{\left(  m-1\right)  \left(  \alpha\left\Vert A\right\Vert
_{1}-\frac{\left\Vert A\right\Vert _{1}^{2}}{mn}\right)  }.
\]
The right-hand side is maximal for $\left\Vert A\right\Vert _{1}=\left(
m+\sqrt{m}\right)  \alpha n/2$ and inequality (\ref{aupb}) follows. If
$\left\Vert A\right\Vert _{1}<n\alpha,$ we see that
\[
\mathcal{E}\left(  A\right)  \leq\sqrt{mtr\left(  AA^{\ast}\right)  }\leq
\sqrt{m\alpha\left\Vert A\right\Vert _{1}}\leq\sqrt{mn}\alpha\leq\alpha
\frac{\left(  m+\sqrt{m}\right)  \sqrt{n}}{2},
\]
completing the proof.
\end{proof}

\textbf{Remarks} \emph{(1)} The bound (\ref{upb}) may be refined using more
sophisticated lower bounds on $\sigma_{1}\left(  A\right)  $. \emph{(2)}
Inequality (\ref{lowb}) and the result of Friedman \cite{Fri04} can be used to
obtain lower bounds for the energy of \textquotedblleft almost
all\textquotedblright\ $d$-regular graphs.


\begin{thebibliography}{99}                                                                                               %


\bibitem {Arn67}L. Arnold, On the asymptotic distribution of the eigenvalues
of random matrices, J. Math. Anal. Appl. 20 (1967) 262--268.

\bibitem {Bal04}R. Balakrishnan, The energy of a graph, Linear Algebra Appl.
387 (2004), 287--295.

\bibitem {Bol98}B. Bollob\'{a}s, \emph{Modern Graph Theory}\textit{,} Graduate
Texts in Mathematics, 184, Springer-Verlag, New York (1998), xiv+394 pp.

\bibitem {CDS80}D. Cvetkovi\'{c}, M. Doob, and H. Sachs, \emph{Spectra of
Graphs,} VEB Deutscher Verlag der Wissenschaften, Berlin, 1980, 368 pp.

\bibitem {Fri04}J. Friedman, A proof of Alon's Second Eigenvalue conjecture, preprint.

\bibitem {FuKo81}Z. F\"{u}redi, J. Koml\'{o}s, The eigenvalues of random
symmetric matrices, Combinatorica 1 (1981), 233--241.

\bibitem {Gut78}I. Gutman, The energy of a graph, Ber. Math.-Stat. Sekt.
Forschungszent. Graz 103 (1978), 1--22.

\bibitem {Gut99}I. Gutman, The energy of a graph: old and new results,
Algebraic Combinatorics and Applications (G\"{o}ssweinstein, 1999), Springer,
Berlin, 2001, pp. 196--211.

\bibitem {HoJo88}R. Horn and C. Johnson, \emph{Matrix Analysis,} Cambridge
University Press, Cambridge, 1985, xiii+561 pp.

\bibitem {KoMo01}J.H. Koolen, V. Moulton, Maximal energy graphs, Adv. Appl.
Math. 26 (2001), 47--52.

\bibitem {KoMo03}J.H. Koolen, V. Moulton, Maximal energy bipartite graphs,
Graphs Combin. 19 (2003), 131--135.

\bibitem {RaTi04}J. Rada, A. Tineo, Upper and lower bounds for the energy of
bipartite graphs, J. Math. Anal. Appl. 289 (2004), 446--455.

\bibitem {Shp06}I. Shparlinski, On the energy of some circulant graphs, Linear
Algebra Appl. 414 (2006), 378--382.

\bibitem {StSt05}D. Stevanovi\'{c}, I. Stankovi\'{c}, Remarks on
hyperenergetic circulant graphs, Linear Algebra Appl. 400 (2005), 345--348.

\bibitem {Wig58}E. Wigner, On the Distribution of the Roots of Certain
Symmetric Matrices, Ann. of Math. 67(1958), 325-328.
\end{thebibliography}
\end{document}